\documentclass[12pt]{article}
\usepackage{amssymb}
\usepackage{amsthm}
\usepackage{amsgen}
\usepackage{amsfonts}
\usepackage{amsbsy}
\usepackage[dvips]{graphicx}

\newtheorem{theorem}{Theorem}[section] 
 
\newtheorem{conjecture}[theorem]{Conjecture} 
\newtheorem{lemma}[theorem]{Lemma} 
\newtheorem{proposition}[theorem]{Proposition} 
 
\newtheorem{example}[theorem]{Example} 
 
\newtheorem{remark}[theorem]{Remark}

\def\pf{{\bf proof}:\ }

\def\qed{$\Box$}

\def\fff{\mathbb{F}}

\def\qqq{\mathbb{Q}}

\def\rrr{\mathbb{R}}

\def\ccc{\mathbb{C}}

\def\zzz{\mathbb{Z}}

\begin{document}

\author{
David Joyner\thanks{Math Dept, USNA, wdj\@@usna.edu, supported by an NSA
Math Sciences grant. }
}
\title{Toric codes over finite fields}
\date{7-26-2003}
\maketitle

\begin{abstract}

In this note, a class of error-correcting codes is
associated to a toric variety associated to a fan defined
over a finite field $\fff_q$, analogous
to the class of Goppa codes associated to a curve.
For such a ``toric code'' satisfying certain
additional conditions, we present an efficient
list decoding algorithm for the dual code.
Many examples are given. For small $q$, many of these 
codes have parameters beating the Gilbert-Varshamov bound.
In fact, using toric codes, we construct a
$(n,k,d)=(49,11,28)$ code over $\fff_8$, which is
better than any other known code listed in 
Brouwer's tables \cite{B} for that $n$ and $k$.
We give upper and conjectural lower bounds on the minimum
distance. The upper bounds are known to be sharp in some cases.
We conclude with a discussion of some decoding methods.

\end{abstract}

\vskip .3in

\section{Introduction}

Some of the constructions discussed here is implemented in 
\cite{MAGMA} and \cite{GAP}, using the {\bf toric}
package\footnote{The GAP version of {\bf toric} requires
\cite{GUAVA}.} \cite{J}.

Let $\fff_q$ denote the finite field with $q$ elements.
A {\bf $q$-ary code} is, for us, a 
a short exact sequence of vector spaces

\begin{equation}
\label{eqn:1}
0 \rightarrow 
C \stackrel{\gamma}{\rightarrow}
\fff^n \stackrel{\theta}{\rightarrow}
\fff^{n-k} \rightarrow
0.
\end{equation}
We identify $C$ with its image under $\gamma$. 
If $\fff^n$ is given the usual standard
vector space basis then the matrix of
$\gamma$ is the generating matrix of $C$
and the matrix of $\theta$ is the parity check matrix 
of $C$.
The integer $n$ is called the
{\bf length} of $C$. The dimension of $C$ is 
$k={\rm dim}_{\fff_q}(C)$. A {\bf codeword} is an
element of $C$. Since the sequence (\ref{eqn:1}) is exact, a
vector $v\in \fff^n$ is a codeword if and only if 
$\theta(v)=0$. 
The {\bf Hamming distance} between two vectors is the number of
coordinates where they differ. This metric stratifies $C$ into
subsets of codewords having the same Hamming distance from 
the origin, or ``weight''. The smallest non-zero weight
which occurs is called the {\bf minimum distance},
denoted $d=d(C)$. 

On one hand, it is known \cite{MS} that
\[
d+k\leq n+1.
\]
On the other hand, it is known \cite{MS} that, as $n$ tends to 
infinity, there are codes
for which $(d/n,k/n)$ lie above the ``Gilbert-Varshamov curve'',
\[
y=1-x\log_q(q-1)-x\log_q(x)-(1-x)\log_q (1-x),
\ \ \ \ \ \ \ \ \ 0\leq x\leq \frac{q-1}{q}.
\]
Many of the toric codes constructed here beat this bound.

\begin{example}
Let $R=\fff_q [x_1,...,x_m]_{\ell}$ denote the
vector space of polynomials in $m$ variables over
$\fff_q$ having degree $\leq \ell$. Let
$S=\fff_q^m=\{p_1,...,p_n\}$ be some indexing,
where $n=q^m$. Let $C$ denote the 
code
\[
C=\{(f(p_1),...,f(p_n))\ |\ f\in R\}.
\]
This code is called the {\bf Reed-Muller code},
denoted ${\cal R}(\ell,m,q)$.

For instance, when $q=2$, $m=5$, $\ell=1$, we have
$n=32$, $k=$dim$_{\fff_2}(C)=6$ and $d=16$.
(This code was used in the 1972 Mariner mission to Mars.)

In general, the length $n$, dimension $k$, and minimum distance $d$
are known exactly for RM codes:
If $\ell=a(q-1)+b$, where $1\leq b\leq q-1$, and if
$0\leq \ell\leq m$, then
$n=q^m$, 
\[
k=
\sum_{i=0}^\ell\sum_{j=0}^{[i/q]}
(-1)^j
\left(
\begin{array}{c}
m\\ j
\end{array}
\right)
\left(
\begin{array}{c}
m-1+i-qj\\ 
m-1
\end{array}
\right),
\]
and $d=(q-b)q^{m-a-1}$.

For this, see for example \cite{MS}, ch 13, \S 3. It should be noted that 
there are easily constructed functions in $R$ of ``minimum weight''
(page 920 in \cite{HCT}).

For example, let $m=2$, $\ell = 5$, and $q$ is ``large''. Then a basis for
$R$ is given by the monomials
$\{x_1^{e_1}x_2^{e_2}\ |\ 0\leq e_1+e_2\leq 5\}$.
We can plot the exponents $(e_1,e_2)$ of these monomials.
It is a triangle with vertices $(0,0)$, $(5,0)$,
$(0,5)$. The number of lattice points in this triangle 
is the dimension of $C$.
\end{example}

All ``toric codes'', defined below, are subcodes of a Reed-Muller code.
They are associated to a complete non-singular $r$-dimensional 
toric variety $X$,
provided with a dense torus $T\hookrightarrow X$, defined
over $\fff_q$, and a $T$-invariant divisor $G$. 
The analogy with the above example is obtained by replacing the triangle
by a convex polytope in $\rrr^r$ and replacing 
all the points in $\fff_q^m$ by a subset of the 
$\fff_q$-rational points in $X$ (often times but not always by
those $\fff_q$-rational
elements in $T$, $T(\fff_q)$). 
The polytope is associated to the divisor $G$
(see \cite{E}, \cite{F}, \cite{JV}, for examples).

We shall focus here on the case of toric surfaces.

\section{Hansen codes}
\label{sec:hansen}

There are three Danish mathematicians working in algebraic
geometry with the same last name, Hansen. 
They are unrelated. In the title of this section,
we are talking about Johan Hansen.

We recall briefly some codes associated to a toric surface,
constructed by J. Hansen \cite{H}.

Let $\fff=\fff_q$ and let 
$\overline{\fff}$ denote a separable algebraic closure. 
Let $L$ be a lattice in $\qqq^2$ generated by $v_1, v_2\in \zzz^2$,
$P$ a polytope in $\qqq^2$,
and $X(P)$ the associated toric surface. Let $P_L=P\cap \zzz^2$.

There is a dense embedding of $GL(1)\times GL(1)$ into $X(P)$
given as follows.
Let $T_L=Hom_{\zzz}(L,GL(1))$ (which is $\cong GL(1)\times GL(1)$
by sending $t=(t_1,t_2)$ to $m_1v_1+m_2v_2\longmapsto
e(\ell)(t)=t_1^{m_1}t_2^{m_2}$)
and let $e(\ell):T_L\rightarrow \overline{\fff}$
be defined by $e(\ell)(t)=t(\ell)$, for $t\in T_L$ and
$\ell\in L$.

Impose an ordering on the set $T_L(\fff )$ (changing the ordering
leads to an equivalent code).
Define the code $C=C_P\subset F^n$ to be the linear code generated
by the vectors
\begin{equation}
\label{eqn:code}
B=\{(e(\ell)(t))_{t\in T_L(F)}\ |\ \ell \in L\cap P_L\},
\end{equation}
where $n=(q-1)^2$. In some special cases, the dimension of $C$ 
and an estimate of its minimum distance can at least be 
conjectured explicitly (see Conjecture \ref{conj:shansen} below).

J. Hansen (\cite{H}, \cite{H3})
gives\footnote{In fact, \cite{H} contains only upper bounds on $d$. 
Explicit computations using the {\bf toric} package suggested
that these upper bounds were attained. This was proven
in \cite{H3}.} 
the minimum distance $d$ of 
such codes in the cases:

(a) $P$ is an isoceles triangle with vertices $(0,0)$,$(a,a)$,$(0,2a)$,

(b) $P$ is an isoceles triangle with vertices $(0,0)$,$(a,0)$,$(0,a)$, or

(c) $P$ is a rectangle with vertices $(0,0)$,$(a,0)$,$(0,b)$,$(a,b)$,

(d) $P$ is a trapazoid with vertices $(0,0)$, $(a,0)$, $(0,b)$, $(a,b+am)$,
where $m>0$.

\noindent
provided $q$ is ``sufficient large''. His precise result is 
recalled below.

\begin{theorem}
\label{thrm:hansen}
Let $a$, $b$ be positive integers. Let $P$ be the polytope 
defined in (a)-(d) above. 

\begin{itemize}
\item[(a)] Assume $q>2a+1$. The code $C=C_P$ has
\[
n=(q-1)^2,\ \ k=(a+1)^2,\ \ d= n-2a(q-1).
\]

\item[(b)] Assume $q>a+1$. The code $C=C_P$ has
\[
n=(q-1)^2,\ \ k=(a+1)(a+2)/2,\ \ d= n-a(q-1).
\]

\item[(c)] Assume $q>{\rm max}(a,b)+1$. The code $C=C_P$ has
\[
n=(q-1)^2,\ \ k=(a+1)(b+1),\ \ d= n-a(q-1)-b(q-1)+ab.
\]

\item[(d)] Assume $q>{\rm max}(a,b,b+am)+1$. The code $C=C_P$ has
\[
n=(q-1)^2,\ \ k=(a+1)(b+1)+m\frac{a(a+1)}{2},\ \ 
\]
\[
d= {\rm min}((q-a-1)(q-b-1),(q-1)(q-b-am-1)).
\]

\end{itemize}
\end{theorem}

The bounds on $q$ are best possible.

We give examples of Hansen's theorem stated above
using MAGMA (or GAP) computations \cite{J}.

\begin{example}

Part (b).

\begin{center}
\begin{tabular}{|c|c|c|c|c|c|}\hline
$q$ & $a$ & $n$ & $k$ & $d$ & using \cite{B}\\ \hline
5 & 2 &  16 & 6 & 8  &Best known\\
5 & 4 &  16 & 13 & 3  & Best possible  \\
5 & 5 &  16 & 15 & 2  & MDS\\

7 & 2 &  36 & 6 & 24  & $d=25$ is best known\\ 
7 & 3 &  36 & 10 & 18  & $d=19$ is best known \\ 

8 & 2 &  49 & 6 & 35  & $d=36$ is best known  \\ 

\hline
\end{tabular}

\end{center}

\end{example}

\section{Construction}
\label{sec:toric}

This section gives a construction which is a little more 
general than Hansen's construction recalled in
\S \ref{sec:hansen}, though it still
falls in the framework of the general class of
codes constructed in \cite{Han2}.

Again, we shall focus on the case of a toric surface,
for simplicity.

Let $M\cong \zzz^2$ be a lattice in $V=\rrr^2$
and let $N\cong \zzz^2$ denote its dual.
Let $\Delta$ be a fan (of rational cones, with
respect to $M$) in $V$ and let $X=X(\Delta)$ denote the
toric variety associated to $\Delta$. 
Let $T$ denote a dense torus in $X$.

Let $\tau_1$, ..., $\tau_s$ denote the edges
of $\Delta$, ordered in a counterclockwise fashion, 
and let $v_i\in N$ denote the 
lattice point closest to the origin in $\tau_i$.
We write $v_{s+1}=v_1$.
The angle between successive vectors $v_i$, $v_{i+1}$ must be
in the range $(0,\pi)$, $1\leq i\leq s$. 
In order for $X$ to be complete, we require that
$\Delta$ contain the cones generated by each succesive pair
$v_i$, $v_{i+1}$, $1\leq i\leq s$. 
In other for $X$ to
be non-singular, we require that each succesive pair
$v_i$, $v_{i+1}$, $1\leq i\leq s$, generate $N$.
(Equivalently, we require that each matrix $(v_i^t,v_{i+1}^t)$
has determinant $\pm 1$.)

Let $P=P_1+...+P_n$ be a positive $1$-cycle on $X$,
where the points $P_i\in X(\fff_q)$ are distinct.
Let $G$ be a $T$-invariant divisor on $X$ which does not ``meet'' $P$,
in the sense that no element of the support of $P$ 
intersects any element in the support of $G$. We
write this as
\[
{\rm supp}(G)\cap {\rm supp}(P)=\emptyset.
\]
Some additional assumptions on $G$ and $P$ shall be made later.
Let 
\[
L(G)=\{0\}\cup \{f\in \fff_q(X)^\times \ |\ div(f)+G\geq 0\}
\]
denote the Riemann-Roch space associated to $G$.

Over $\ccc$, according to \cite{F}, \S 3.4, there is a polytope 
$P_G$ in $V$ such that $L(G)$ is spanned by the
monomials $x^a$ (in multi-index notation), for $a\in P_G\cap N$.
Over $\fff_q$, we assume that the same result holds provided that,
for all $a=(a_1,...,a_r)\in P_G\cap N$, we have 
$|a_i|<q$, for all $1\leq i\leq r$. {\bf In other words, we assume $q$
is ``sufficiently large'', in this sense (depending on $G$).}

Let $C_L=C_L(P,G)$ denote the code define by
\begin{equation}
\label{eqn:code2}
C_L=\{(f(P_1),...,f(P_n))\ |\ f\in L(G)\}.
\end{equation}
This is the {\bf Goppa code associated to $X$, $P$,
and $G$}. The dual code is denoted
\begin{equation}
\label{eqn:dualcode}
C=C_L^\perp = \{(c_1,...,c_n)\in \fff_q^n\ |\ \sum_{i=1}^n c_if(P_i)=0,
\forall f\in L(G)\}.
\end{equation}

It is reasonable to ask under what conditions (if any) is the map
\begin{equation}
\label{eqn:evalmap}
\begin{array}{ccc}
L(G) & \rightarrow &C_L\\
f & \longmapsto & (f(P_1),...,f(P_n))
\end{array}
\end{equation}
an injection? This question was investigated by S. Hansen in \cite{Han3}
and \cite{Han1} for varieties, \S 5.3, and \cite{H3} 
in the case of certain toric surfaces.

\begin{lemma} (J. Hansen)
\label{lemma:eval}
If $X=X(\Delta)$ is a non-singular toric surface and $q$ is 
sufficiently large (i.e., $q>q_0(\Delta)$, for some $q_0(\Delta)>1$)
then (\ref{eqn:evalmap}) is injective.
\end{lemma}

This was proven in \S 2.3 in \cite{H3} for certain cases
(where $q_0$ was explicitly determined)
but the technique applies more generally to yield the
above result.

Let $\Delta$ be a fan in a lattice $L\cong \zzz^r$. 
Let $\tau_i$ and $v_i$ be as above $1\leq i\leq s$.
Let $D_i$ denote the Weil divisor 
\begin{equation}
\label{eqn:divisor}
D_i={\rm Hom}(\tau_i^\perp\cap L^\perp,\ccc^\times),
\end{equation}
which can be regarded as the Zariski closure of an orbit of $T$. 

\begin{example}
\label{example:fan7}

Let $\Delta$ be the fan generated by cones
whose edges are defined by
\[
v_1=5e_1-e_2,\ \ v_2=-e_1+5e_2,\ \ v_3=-e_1-e_2.
\]
Let $X$ be the toric variety associated to $\Delta$.
In this example, the divisor $G=5D_3$ on $X/\fff_8$ yields a 
toric code which one finds (thanks to the {\bf toric}
package) has the parameters $(49,11,28)_8$. 
This is a new code and beats the previous best known code
(for given $n$, $k$, $q$)
having parameters $(49,11,27)_8$ \cite{B}.

\end{example}

\begin{example}
\label{example:fan1}

Let $\Delta$ be the fan generated by cones
whose edges are defined by
\[
v_1=2e_1-e_2,\ \ v_2=-e_1+2e_2,\ \ v_3=-e_1-e_2.
\]
Let $X$ be the toric variety associated to $\Delta$.

\vskip .3in

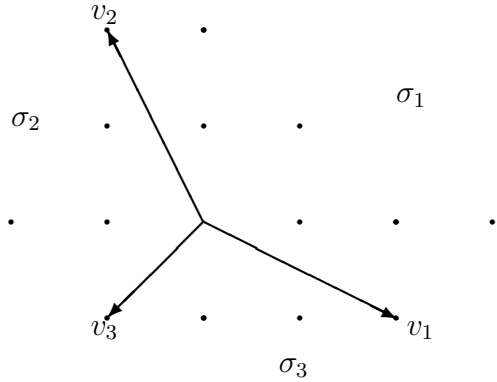
\begin{center}
\setlength{\unitlength}{.02cm}
\begin{picture}(256.00,150.00)(-120.00,0.00)
\thicklines
\put(0.00,64.00){\circle*{4}} 
\put(0.00,128.00){\circle*{4}} 

\put(-64.00,0.00){\circle*{4}} 
\put(0.00,-64.00){\circle*{4}} 

\put(-64.00,-64.00){\circle*{4}} 
\put(-64.00,64.00){\circle*{4}} 
\put(-64.00,128.0){\circle*{4}} 

\put(64.00,0.00){\circle*{4}} 
\put(128.00,0.00){\circle*{4}} 
\put(192.00,0.00){\circle*{4}} 

\put(0.00,0.00){\vector(-1,2){64.00}} 
\put(0.00,0.00){\vector(-1,-1){64.00}} 

\put(64.00,64.00){\circle*{4}} 
\put(128.00,0.00){\circle*{4}} 
\put(0.00,128.00){\circle*{4}} 
\put(-128.00,0.00){\circle*{4}} 

\put(128.00,80.00){$\sigma_1$} 
\put(-128.00,64.00){$\sigma_2$} 
\put(50.00,-100.00){$\sigma_3$} 

\put(135.00,-75.00){$v_1$} 
\put(-75.00,135.00){$v_2$} 
\put(-75.00,-75.00){$v_3$} 

\put(0.00,-64.00){\circle*{4}} 
\put(64.00,-64.00){\circle*{4}} 
\put(128.00,-64.00){\circle*{4}} 

\put(0.00,0.00){\vector(2,-1){129.00}} 
\end{picture}
\end{center}

\vskip 1in

\noindent
Note $v_1,v_2$ generate a cone $\sigma_1$ in the 
fan $\Delta$ yet do not generate the lattice $\zzz^2$.
By the criterion on page 29 of \cite{F}, $X$ is singular.

In the notation of (\ref{eqn:divisor}), 
the divisor $G=d_1D_1+d_2D_2+d_3D_3$
is a Cartier divisor if and only if
$d_1\equiv d_2\equiv d_3\ ({\rm mod}\ 3)$ (this is
a consequence of an exercise on page 62 of \cite{F}).
A Cartier divisor is very ample if and only if
$d_1+d_2+d_3> 0$ (this follows from the
criterion in the proof of the proposition on page
68 in \cite{F}). Let 
\[
\begin{array}{c}
P_G=\{(x,y)\ |\ \langle (x,y),v_i\rangle \geq -d_i, \ \forall i\}\\
=\{(x,y)\ |\ 2x-y \geq -d_1, -x+2y \geq -d_2, -x-y \geq -d_3 
\}
\end{array}
\]
denote the polytope associated to $G$. 

We tabulate parameters of toric codes
$C_L(P,G)$, where $P$ is the sum of all the
$\fff_q$-valued points in a ``dense'' torus $T$ in $X$.
These results were obtained using the {\bf toric}
package. Note that not all of these
divisors $G$ are Cartier. 

\begin{center}
\begin{tabular}{|c|c|c|c|c|c|c|c|}\hline
$q$ & $d_1$ & $d_2$ & $d_3$ & $n$ & $k$ & $d$ & Using \cite{B} \\ \hline
 5   &  0     &  0     &  3     &  16   &  4   &  10  & $d=11$ is best possible \\
 5   &  1     &  1     &  2     &  16   &  5   &  8  & $d=9$ is best known \\
 5   &  1     &  1     &  3     &  16   &  7   &  6  &$d=7$ is best known  \\
 5   &  1     &  1     &  4     &  16   &  10   &  4  & $d=5$ is best possible \\
 5   &  2     &  3     &  1     &  16   &  9   &  6   &Best possible   \\
 5   &  2     &  3     &  2     &  16   &  12   &  3  & $d=4$ is best possible \\
 5   &  2     &  3     &  3     &  16   &  14   &  2  & Best possible   \\
 5   &  2     &  3     &  3     &  16   &  15   &  2   & MDS\\

 7   &  1     &  2     &  0     &  36   &   3  &  30 & Best possible   \\
 7   &  1     &  2     &  1     &  36   &   5  &  24 & $d=26$ is best known  \\
 7   &  1     &  2     &  2     &  36   &   7  &  21 & $d=22$ is best known  \\
 7   &  1     &  2     &  3     &  36   &   9  &  20 & Best known  \\
 7   &  2     &  3     &  4     &  36   &   18  &  9 &  $d=12$ is best known  \\
 7   &  3     &  0     &  0     &  36   &   4  &  27 & $d=28$ is best known  \\
8   &  0  &  0  &  3  &  49  &  4  &  40  & Best known \\
8   &  0  &  1  &  2  &  49  &  3  &  42  & Best possible \\
8   &  0  &  2  &  3  &  49  &  7  &  33  & $d=35$ is best known \\
8   &  0  &  3  &  3  &  49  &  10  &  28  & Best known \\
8   &  0  &  3  &  4  &  49  &  12  &  26  & Best known \\

\hline
\end{tabular}

\end{center}

\end{example}

\begin{example}
\label{example:fan5}

Let for now $\Delta$ be the fan generated by 
\[
v_1=e_1,\ \ v_2=-e_1+e_2,\ \ v_3=-e_1-e_2.
\]
Let $X$ be the toric variety associated to $\Delta$.
This is singular.
The divisor $G=d_1D_1+d_2D_2+d_3D_3$ with
$d_1=0$, $d_2=0$, $d_3=a$ gives case (a) of
Hansen's code.

To create a smooth toric variety yielding Hansen's code,
let instead $\Delta$ be the refined fan generated by 
\[
v_1=e_1,\ \ v_2=-e_1+e_2,\ \ v_3=-e_1-e_2,\ \ v_4=-e_2.
\]
Let $X$ be the toric variety associated to $\Delta$.
This is now non-singular.
The divisor $G=d_1D_1+d_2D_2+d_3D_3+d_4D_4$ with
$d_1=0$, $d_2=0$, $d_3=a$, $d_4=2a$ (for example)
gives case (a) of Hansen's code.
\end{example}

\begin{example}
\label{example:fan4}

Let $\Delta$ be the fan generated by 
\[
v_1=e_1,\ \ v_2=e_2,\ \ v_3=-e_1-e_2.
\]
Let $X$ be the toric variety associated to $\Delta$.
This is non-singular.
The divisor $G=d_1D_1+d_2D_2+d_3D_3$ with
$d_1=0$, $d_2=0$, $d_3=a$ gives case (b) of
Hansen's code.
\end{example}

\begin{example}
\label{example:fan3}

Let $\Delta$ be the fan generated by 
\[
v_1=e_1,\ \ v_2=e_2,\ \ v_3=-e_1,\ \ v_4=-e_2.
\]
Let $X$ be the toric variety associated to $\Delta$.
This is non-singular.
The divisor $G=d_1D_1+d_2D_2+d_3D_3+d_4D_4$ with
$d_1=0$, $d_2=0$, $d_3=a$, $d_4=b$ gives case (c) of
Hansen's code.

\end{example}

\begin{example}
\label{example:fan2}

Let $\Delta$ be the fan generated by 
\[
v_1=e_1,\ \ v_2=-e_1+me_2\ \ (m>1)\ {\rm fixed}),\ \ v_3=-e_2.
\]
Let $X$ be the toric variety associated to $\Delta$.
This surface is, except for small $m$, singular.

In the notation above, the divisor $G=d_1D_1+d_2D_2+d_3D_3$
is a Cartier divisor if and only if
$d_1 + d_2\equiv 0\ ({\rm mod}\ m)$ (this is
an exercise on page 64 of \cite{F}).

Let $\Delta'$ be the refined fan generated by 
\[
v_1=e_1,\ \ v_2=e_2,\ \ v_3=-e_1+me_2\ \ (m>1)\ {\rm fixed}),
\ \ v_4=-e_2.
\]
Let $X'$ be the toric variety associated to $\Delta'$.
This surface is non-singular.
The divisor $G'=d'_1D'_1+d'_2D'_2+d'_3D'_3+d'_4D_4$
(using $D'_i$ to denote the divisor in
(\ref{eqn:divisor}) defined on $X'$)
is very ample if and only if
$d'_2+d'_4\geq 0$, $d'_1+d'_3\geq md'_1$ (this is an exercise 
on page 70 of \cite{F}).
Taking 
\[
d_1'=d_2'=0, \ \ d_3'=a,\ \ d_4'=b,
\]
and swapping the roles of $x$ and $y$, gives the Hansen code
of type (d) constructed above.

Let 
\[
\begin{array}{c}
P_G=\{(x,y)\ |\ \langle (x,y),v_i\rangle \geq -d_i, \ \forall i\}\\
=\{(x,y)\ |\ x \geq -d_1, -x+my \geq -d_2, -y \geq -d_3 
\}
\end{array}
\]
denote the polytope associated to the Weil
divisor $G=d_1D_1+d_2D_2+d_3D_3$.

The {\bf toric} package gives the following results:

$m=3$

\begin{center}
\begin{tabular}{|c|c|c|c|c|c|c|c|}\hline
$q$ & $d_1$ & $d_2$ & $d_3$ & $n$ & $k$ & $d$ & Using \cite{B} \\ \hline
5    & 0  &  2 & 2  & 16  & 11  & 3 & $d=4$ is best known  \\
5    & 0  &  2 & 3  & 16  & 15  & 2 & MDS  \\
5    & 0  &  4 & 2  & 16  & 14  & 2 & Best possible  \\
5    & 1  &  0 & 0  & 16  & 2  & 12 & $d=13$ is best possible  \\
7    & 0  &  1 & 0  & 36  & 2  & 30 & $d=31$ is best possible  \\
8    & 1  &  0 & 0  & 49  & 2  & 42 & $d=43$ is best possible  \\
9    &  0 &  1 & 0  & 64  & 2  & 56 & $d=57$ is best possible  \\
\hline
\end{tabular}
\end{center}

$m=5$ 

\begin{center}
\begin{tabular}{|c|c|c|c|c|c|c|c|}\hline
$q$ & $d_1$ & $d_2$ & $d_3$ & $n$ & $k$ & $d$ & Using \cite{B} \\ \hline
5    & 0  &  0 & 3  & 16  & 13 & 2 & $d=3$ is best possible  \\
5    & 3  &  3 & 2  & 16  & 14 & 2 & Best possible  \\
7    & 1  &  3 & 4  & 36  & 29  & 3 & $d=5$ is best known  \\
8    & 4  &  4 & 4  & 49  & 39  & 3 & $d=6$ is best known  \\
\hline
\end{tabular}
\end{center}

$m=10$ 

\begin{center}
\begin{tabular}{|c|c|c|c|c|c|c|c|}\hline
$q$ & $d_1$ & $d_2$ & $d_3$ & $n$ & $k$ & $d$ & Using \cite{B} \\ \hline
7    & 5  &  7 & 4  & 36  & 33  & 2 & Best possible \\
8    & 5  &  9 & 4  & 49  & 40  & 3 & $d=6$ is best known  \\
9    & 5  &  9 & 4  & 64  & 45  & 4 & $d=12$ is best known \\

\hline
\end{tabular}
\end{center}

Indeed, the parameters arising in these case are very similar
to each other.
\end{example}

\begin{example}
\label{example:fan6}

Let $\Delta$ be the fan generated by 
\[
v_1=2e_1-e_2,\ \ v_2=-e_1+e_2,\ \ v_3=-e_1.
\]
Let $X$ be the toric variety associated to $\Delta$.
This is non-singular.
All the divisors $G=d_1D_1+d_2D_2+d_3D_3$
are Cartier.
The divisor $G=d_1D_1+d_2D_2+d_3D_3$
is very ample if and only if $d_1+d_2+d_3>0$.

Let 
\[
\begin{array}{c}
P_G=\{(x,y)\ |\ \langle (x,y),v_i\rangle \geq -d_i, \ \forall i\}\\
=\{(x,y)\ |\ 2x-y \geq -d_1, -x+y \geq -d_2, -x \geq -d_3 
\}
\end{array}
\]
denote the polytope associated to $G$. 
If each $d_i>0$ and if $d_3>d_1+d_2$ then the volume of this polytope 
is given by
\[
{\rm vol}(P_G)=\frac{1}{2}(2d_1+d_3)(d_3-d_1-d_2).
\]

We tabulate parameters of toric codes
$C_L(P,G)$, where $P$ is the sum of all the
$\fff_q$-valued points in a ``dense'' torus $T$ in $X$.
These parameters were obtained using the {\bf toric} package. 

\begin{center}
\begin{tabular}{|c|c|c|c|c|c|c|c|}\hline
$q$ & $d_1$ & $d_2$ & $d_3$ & $n$ & $k$ & $d$ & Using \cite{B} \\ \hline
5    &  0     &  0     &  1     &  16   &  3   &  12 & Best possible   \\
5    &  0     &  0     &  2     &  16   &  6   &  8  & $d=9$ is best known \\
5    &  0     &  0     &  3     &  16   &  10   &  4  & $d=5$ is best known \\
5    &  0     &  0     &  4     &  16   &  13   &  3 & Best possible   \\
7    &  0     &  0     &  1     &  36   &  3   &  30  & Best possible   \\
7    &  0     &  0     &  2     &  36   &  6   &  24  & Best known \\
7    &  0     &  0     &  3     &  36   &  10   &  18  & Best known \\
7    &  0     &  0     &  4     &  36   &  15   &  12  & $d=14$ is best known \\
7    &  4     &  1     &  1     &  36   &  26   &  5  & $d=6$ is best known \\
7    &  4     &  1     &  2     &  36   &  30   &  4 & $d=5$ is best possible    \\
7    &  4     &  1     &  3     &  36   &  33   &  3 & Best possible   \\
7    &  4     &  1     &  4     &  36   &  35   &  2  & MDS \\
8    &  0     &  0     &  1     &  49   &  3   &   42  & Best possible \\
8    &  0     &  0     &  2     &  49   &  6   &   35  & $d=36$ is best known   \\
8    &  0     &  0     &  3     &  49   &  10   &  28  & Best known \\
8    &  0     &  0     &  4     &  49   &  15   &  21  & $d=23$ is best known \\
8    &  0     &  4     &  3     &  49   &  34   &  6  & $d=10$ is best known \\
8    &  0     &  4     &  4     &  49   &  39   &  5 & $d=6$ is best known   \\
8    &  2     &  4     &  4     &  49   &  46   &  3 & best possible   \\
8    &  3     &  4     &  4     &  49   &  48   &  2 & MDS  \\
8    &  4     &  1     &  4     &  49   &  43   &  4 & $d=5$ is best possible   \\
9    &  0     &  0     &  1     &  64   &   3  &   56  & Best possible  \\
9    &  0     &  0     &  2     &  64   &   6  &   48  & $d=49$ is best known  \\
9    &  0     &  0     &  3     &  64   &   10  &  40  & $d=41$ is best known \\

\hline
\end{tabular}

\end{center}

\end{example}

\section{Estimates on the parameters}
\label{sec:estimate}

We shall now estimate the parameters $n,k,d$ for 
$C_L=C_L(G,P,X)$, conjecturally. Assume $X$ is a complete toric
variety of dimension $r$. For many such varieties,
the conjectured estimates below were verified by computer 
using the {\bf toric} package.

Let $G$ be a $T$-invariant Cartier divisor on $X=X(\Delta)$,
let $P_G$ be the polytope associated to $G$ and let $M$ be 
the lattice as in \cite{F}.

\begin{conjecture}
\label{conj:shansen}
Let $C_L=C_L(G,P,X)$ be as in (\ref{eqn:code2}). Assume 
\begin{itemize}
\item
$X$ is a non-singular,
projective toric variety of dimension $r$, 
\item
$n$ is so large that there is an integer $N>1$ such that
\[
2N\cdot {\rm vol}(P_G)\leq n\leq 2N^2\cdot {\rm vol}(P_G).
\]
\end{itemize}
If $q$ is ``large'' then any 
$f\in L(G)=H^0(X,{\cal O}(G))$ has no more than $n$ zeros
in $X(\fff_q)$.
Consequently, 

\[
d\geq n- 2N\cdot {\rm vol}(P_G).
\]

\end{conjecture}

Here ``large'' make depend on $X$, $P$ and $G$ but not
on $f$.

Let $G$ be a $T$-invariant Cartier divisor on $X=X(\Delta)$
and let
\[
\psi_G:|\Delta|\rightarrow \rrr
\]
be the linear function associated to $G$ as in \cite{F},
\S 3.4 (this is called a {\bf virtual support
function} in \cite{E}, Definition VII.4.1).
If $\psi_G$ is strictly convex and $X$ is non-singular
then $G$ is very ample
and ${\cal O}(G)$ is generated by its global sections
(\cite{F}, \S 3.4, p. 70)\footnote{
Demazure has shown that on a complete, non-singular toric variety
a $T$-invariant divisor is ample if and only if it is very ample
(see \cite{F}, page 71).
}.
Conversely, if ${\cal O}(G)$ is generated by its global sections then
$\psi_G$ is convex and 
\[
\psi_G(v)={\rm min}_{u\in P_G\cap M}\ \langle u,v\rangle,
\]
where $P_G$ and $M$ are as above.

\begin{conjecture}
\label{conj:main}
Let $C_L=C_L(G,P,X)$ be as in (\ref{eqn:code2}). Assume 
\begin{itemize}
\item
$X$ is a non-singular,
projective toric variety of dimension $r$, 
\item
$\psi_G$ is strictly convex,
\item
and  
\[
{\rm deg}(P)>{\rm deg}(G^r).
\]
\end{itemize}
If $q$ is ``large'' then any 
$f\in L(G)=H^0(X,{\cal O}(G))$ has no more than $n={\rm deg}(P)$ zeros
in $X(\fff_q)$.
Consequently, 

\[
k\geq {\rm dim}\, H^0(X,{\cal O}(G))=|P_G\cap M|,
\]
and
\[
d\geq n-r!\cdot |P_G\cap M|.
\]
Moreover, if $n>r! |P_G\cap M|$ then 
${\rm dim}\, H^0(X,{\cal O}(G))=|P_G\cap M|$.
\end{conjecture}

\begin{remark}
This is false if $X$ is singular.
\end{remark}

We conclude with an upper bound on the minumum distance.
Let $M^+=\{m\in M\ |\ m_i\geq 0\}$.

\begin{theorem}
Let $C=C_L(G,P,X)$ be as in (\ref{eqn:code2}). Let $H\subset P_G\cap M^+$
be a line $\{(e_1,...,e_{j-1},x,e_{j+1},...,e_r)\ |\ 0\leq x\leq h\}$,
for some integer $h>0$.
Then $d\leq n-h$.

\end{theorem}

\begin{remark}
(1) In particular, let $r=2$ and let $h>0$ denote the 
largest integer such that $\{(0,0),...,(0,h)\}\subset P_G$.
Then $d\leq n-h$.

(2) This generalizes the results of J. Hansen stated above,
which are known to be sharp.

\end{remark}

\pf
Let $h={\rm max}_i(d_i)$. Let $f\in L(G)$ denote the function
\[
(x_j-a_1)...(x_j-a_h)\prod_{i\not= j}x_i^{e_i}\in L(G),
\]
where $a_i\in \fff_q^\times$ are distinct.
This vanishes at $h$ points, so $d\leq n-h$.
\qed

\section{Decoding algorithm}
\label{sec:decode}

One reasonable decoding algorithm is an algorithm for
decoding more general multi-dimensional cyclic
codes. We refer to J. Little's papers \cite{L1}, \cite{L2} 
and the references cited there. 

In this section, we discuss a much more restrictive algorithm.
The algorithm here will decode $C=C_L(P,G)^\perp$ when
$G$ is a ``large'' divisor and the number of errors is ``small''.
What it lacks in optimality it makes up in simplicity.

Let $X=X(\Delta)$ be a toric variety associated to a fan $\Delta$
with edges $\tau_1$, ..., $\tau_s$. Let $D_i$ denote the $T$-Weil
divisors as above.

Let $C$ be as in (\ref{eqn:dualcode}).
Let $r=c+e$ be a received vector, $c\in C$ a codeword
and $e$ the error vector of smallest weight. Thus all
$f\in L(G)$ satisfy $\sum_{i=1}^n c_if(P_i)=0$. Let
\[
I=\{i\ |\ e_i\not= 0\}.
\]

We shall follow \cite{S}, pages 217-218, to (a) create an
error locator function for $r$ and $C$, and (b) solve for 
$e_i$, $i\in I$. The method generally only gives a list
of several vectors, one of which is the true error
vector $e$ - it only yields one vector (namely, $e$)
under a certain technical condition (see
Proposition \ref{prop:errors} for details).

Let $G'$ be another divisor of $X$. We shall make some
assumptions on $G$ and $G'$ as we go along. Let
\[
L(G')={\rm span}\{f_1,...,f_\ell\},
\]
\[
L(G-G')={\rm span}\{g_1,...,g_k\},
\]
\[
L(G)={\rm span}\{h_1,...,h_m\}.
\]

We assume 
\begin{itemize}
\item[(A)]
that these spaces $L(G),L(G'),L(G-G')$ are non-zero,
\item[(B)]
that there are $T$-invariant divisors $D_i$ for which 
(a) each $P_i$ is contained in a $D_i$, $i\in I$,
(b) $P_i\notin D_j$, for $i\not= j$ and $i,j\in I$,
(c) $L(G'-\sum_{i\in I}D_i)$ is non-zero,
(d) for any fixed $i_0\in I$,
$L(G-G'-\sum_{i\in I}D_i)$ is non-empty and {\it properly}
contained in 
$L(G-G'-\sum_{i\in I,i\not= i_0}D_i)$.
\item[(C)]
$d> {\rm deg}(G')$.

\end{itemize}
Assumption (B) above is very strong and forces us
to assume that $|I|\leq s$. Roughly speaking, $G$ is a ``large'' 
divisor and the number of errors is ``small''.

Let
\[
[r,\phi ]=\sum_{i=1}^n r_i\phi (P_i),
\]
and let $x_1=a_1$, ..., $x_\ell=a_\ell$ be a non-trivial solution 
to the system of $k$ simultaneous equations
\begin{equation}
\label{eqn:coeffs1}
\sum_{j=1}^\ell [r,f_jg_i]x_j=0,
\ \ \ \ \ \ \ 1\leq i\leq k.
\end{equation} 
(We show such a solution exists below.)
Let $f=a_1f_1+...+a_\ell f_\ell$.

\begin{proposition}
\label{prop:errorlocator}
Under the assumptions above, 
$f$ is an error locator function. In other words,
$f(P_i)=0$ for all $i\in I$.
\end{proposition}

\pf
We have $L(G'-\sum_{i\in I}D_i)$ is non-zero, by assumption.
Choose a non-zero $z\in L(G'-\sum_{i\in I}D_i)\subset L(G')$
and write $z=\gamma_1f_1+..+\gamma_\ell f_\ell$, for some
$\gamma_i\in \fff_q$. Then $zg_j\in L(G)$, for all
$1\leq j\leq k$. If $r=(r_1,...,r_\ell)$ then 
\[
[r,zg_j]=\sum_{i=1}^\ell [r,f_ig_j]\gamma_i.
\]
On the other hand, $[r,zg_j]=[c+e,zg_j]=[e,zg_j]$, so
\[
[r,zg_j]=[e,zg_j]
= \sum_{v=1}^n e_v f_i(P_v)g_j(P_v).
\]
But $e_v=0$ for $v\notin I$ and $z(P_i)=0$ for
$v\in I$, so $[r,zg_j]=0$. This implies that
$(\gamma_1,...,\gamma_\ell)$ is a non-trivial solution to 
(\ref{eqn:coeffs1}).

Now suppose that there is an $i_0\in I$ for which 
$f(P_{i_0})\not= 0$. By assumption (B,d),
there is an $h\in L(G-G')$ such that
$h(P_{i_0})\not= 0$ and $h(P_i)=0$ for all
other $i\in I$. Then
\[
[r,fh]=\sum_{v=1}^n e_v f(P_v)h(P_v)
=e_{i_0}f(P_{i_0})h(P_{i_0})\not= 0.
\]
On the other hand, $h$ is a linear combination 
of the $g_j$'s. So this contradicts the 
equation $[r,zg_j]=0$ we obtained earlier.
\qed

Next, we must determine the $e_i$, $i\in I$.
Let 
\[
N(f)=\{i \ |\ 1\leq i\leq n, f(P_i)=0\}.
\]
By the proposition, $I\subset N(f)$. 

\begin{proposition}
\label{prop:errors}
The $b_i=e_i$, for $i\in N(f)$, solve 
\[
\sum_{i\in N(f)} b_i h_j(P_i)=[r,h_j],\ \ \ \ \ 1\leq j\leq m.
\]
If $|N(f)|\leq {\rm deg}(G')$ then this solution is unique.
\end{proposition}

\begin{remark}
Naturally, it would be interesting to know if
the hypothesis $|N(f)|\leq {\rm deg}(G')$ can be removed.
Does the
$T$-invariance of the $D_i$ and the fact that each
$D_i$ is (the Zariski closure of) an orbit help?

A bound on $|N(f)|$ is known for certain $\Delta$ - 
see \S 2.3 in \cite{H3}. A conjectural bound is given in the
previous section.
\end{remark}

\pf 
Since $h_j\in L(G)$ and (by the previous proposition)
$e_v=0$ for all $v\notin N(f)$, we have
\[
[r,h_j]=[e,h_j]=\sum_{v=1}^n e_v h_j(P_v)
=\sum_{v\in N(f)} e_v h_j(P_v).
\]
Thus the equation displayed in Proposition
\ref{prop:errors} has the $e_i$'s, $i\in N(f)$, as a solution.

Now we show that this solution is unique. Suppose
$b_i$, $i\in N(f)$, is another solution. Define
$b_j=0$ for $j\in \{1,2,...,n\}-N(f)$. Then
for $1\leq j\leq m$,
\[
[b,h_j]=\sum_{v=1}^n b_v h_j(P_v)
=\sum_{v\in N(f)} b_v h_j(P_v)=[r,h_j]=[e,h_j].
\]
As the $h_j$'s form a basis of $L(G)$, this implies
$b-e$ is a codeword. The weight of this codeword
satisfies 
\[
wt(b-e)\leq |N(f)|\leq deg(G')< d,
\]
by our hypothesis. This forces $b_i=e_i$,
so the solution is unique. 
\qed

The above Proposition tells us how to find the error vector, given 
$f$, $L(G)$, and the $P_i$'s, as desired.

We illustrate this algorithm with an example.

\begin{example}

Here is an example where we decode 
a received word in a toric code with $\leq 3$ errors.

Let $\Delta$ be the fan as in Example \ref{example:fan1}.
Let $X$ be the toric variety associated to $\Delta$.

Let 
\[
\begin{array}{c}
P_D=\{(x,y)\ |\ \langle (x,y),v_i\rangle \geq -d_i, \ \forall i\}\\
=\{(x,y)\ |\ 2x-y \geq -d_1, -x+2y \geq -d_2, -x-y \geq -d_3 
\}
\end{array}
\]
denote the polytope associated to the Weil
divisor $D=d_1D_1+d_2D_2+d_3D_3$,
where $D_i$ is as above. 

Let 
\[
G=10D_3,\ \ \ 
D=D_1+D_2+D_3,\ \ \ G'=2D.
\]
Then $P_G$ is a triangle in the plane 
with vertices at $(0,0)$, $(10/3,20/3)$, and $(20/3,10/3)$. 
It's area is $\frac{1}{2}bh=50/3$.

Let $T\subset X$ denote the dense torus of $X$.
In this example, the patch $U_{\sigma_1}$ is an affine variety with 
coordinates $x_1,x_2,x_3$ given by $x_1^3-x_2x_3=0$. 
The torus embedding
$T\hookrightarrow U_{\sigma_1}$ is given by
sending $(t_1,t_2)$ to 
$(x_1,x_2,x_3)=(t_1t_2,t_1t_2^2,t_1^2t_2)$.
The patch $U_{\sigma_2}$ is an affine variety with 
coordinates $y_1,y_2,y_3$ given by $y_2^2-y_1y_3=0$. 
The torus embedding
$T\hookrightarrow U_{\sigma_2}$ is given by
sending $(t_1,t_2)$ to 
$(y_1,y_2,y_3)=(t_1^{-2}t_2^{-1},t_1^{-1},t_1^{-1}t_2)$.
The patch $U_{\sigma_3}$ is an affine variety with 
coordinates $z_1,z_2,z_3$ given by $z_2^2-z_1z_3=0$. 
The torus embedding
$T\hookrightarrow U_{\sigma_3}$ is given by
sending $(t_1,t_2)$ to 
$(x_1,x_2,x_3)=(t_1^{-1}t_2^{-2},t_2^{-1},t_1t_2^{-1})$.

In the local coordinates of $U_{\sigma_1}$, the
space $L(G)$ has as a basis,
\[
\begin{array}{c}
\{h_i\}=
\{ 1, x_1 x_2, x_1^2 x_2, x_1 x_2^2, x_1^2 x_2^2, 
x_1^3 x_2^2, x_1^4 x_2^2, x_1^2 x_2^3, x_1^3 x_2^3,\\ 
x_1^4 x_2^3, x_1^5 x_2^3, x_1^6 x_2^3, x_1^2 x_2^4, 
x_1^3 x_2^4, x_1^4 x_2^4, x_1^5 x_2^4, x_1^6 x_2^4, \\
x_1^3 x_2^5, x_1^4 x_2^5, x_1^5 x_2^5, x_1^3 x_2^6, 
x_1^4 x_2^6\}\, .
\end{array}
\]
In particular, it is 22-dimensional.

The polytope $P_{G'}$ has vertices at
$(-2,-2)$, $(2,0)$, $(0,2)$.
In the local coordinates of $U_{\sigma_1}$, the
space $L(G')$ has as a basis,
\[
\{f_i\}=
\{ x_1^{-2}x_2^{-2}, x_1^{-1}, x_2^{-1}, x_1^{-1}x_2^{-1},  1, x_1, x_1^2, 
x_1x_2, x_2, x_2^2\}.
\]
In particular, it is 10-dimensional.

The polytope $P_{G-G'}$ has vertices at
$(2,2)$, $(10/3,14/3)$, $(14/3,10/3)$.
In the local coordinates of $U_{\sigma_1}$, the
space $L(G-G')$ has as a basis,
\[
\{g_i\}=
\{x_1^2 x_2^2, x_1^3 x_2^3, x_1^4 x_2^3, x_1^3 x_2^4, x_1^4 x_2^4\}.
\]
$L(G-G')$ is 5-dimensional.

Choose distinct points $P_1,...,P_n\in X(\fff_q)$ and let
\[
C=\{(c_1,...,c_n)\in \fff_q^n\ |\ \sum_{i=1}^n c_i f(P_i)=0,
\forall f\in L(G)\}.
\]

Let $P_i\in D_i$ and $P_i\notin D_j$, with
$1\leq i\not= j\leq 3$. Assume $P_4,...,P_n\notin D_i$, $i=1,2,3$.
(For example, let $\{P_4,...,P_n\}\subset T(\fff_q)$.)
We must choose $n>50$, assuming 
Conjecture \ref{conj:main} below, to insure condition
(C) above is satisfied\footnote{For 
example, we know
\[
\begin{array}{c}
|X(\fff_2)|=7,\ \ |X(\fff_3)|=13,\ \ |X(\fff_4)|=21,\ \\  
|X(\fff_5)|=31,\ \ |X(\fff_7)|=57,\ \  |X(\fff_8)|=73,...
\end{array}
\]
So, we may take $q\geq 7$, at least conjecturally. See 
also Hansen's Lemma \ref{lemma:eval}.
}.

Since $L(G-3D)$ is 1-dimensional, hence non-zero,
the hypotheses and assumptions above are all
satisfied. 

Suppose $r=(1,1,1,0,...,0)\in\fff_q^n$ is the
received vector (it has $3$ errors in the first $3$ positions).
In this case,
\[
[r,f]=f(P_1)+f(P_2)+f(P_3).
\]
Let $a_1$, $a_2$, ..., $a_{10}$ be a solution to
\[
\begin{array}{c}
a_1[r,f_1g_1]+a_2[r,f_2g_1]+...+a_{10}[r,f_{10}g_1]=0,\\
a_1[r,f_1g_2]+a_2[r,f_2g_2]+...+a_{10}[r,f_{10}g_2]=0,\\
\vdots\\
a_1[r,f_1g_{5}]+a_2[r,f_2g_{5}]+..+a_{10}[r,f_{10}g_{5}]=0.
\end{array}
\]
This determines
\[
f=\sum_{i=1}^{10} a_jf_j,
\]
and hence $N(f)$.
The equations
\[
\sum_{i\in N(f)} b_i h_j(P_i)=[r,h_j],\ \ \ \ \ 1\leq j\leq 22,
\]
clearly have a solution ($b_1=b_2=b_3=1$, the other $b_i=0$).

\end{example}

{\bf Acknowledgements}:
I thank Soren Hansen, Johan Hansen, and Amin Shokrollahi for 
answering many questions and helpful
correspondence regarding their research.


\begin{thebibliography}{99}  


\bibitem[B]{B} W. Brouwer's tables on
Bounds on the minimum distance of linear codes,
\newline
\verb+http://www.win.tue.nl/~aeb/voorlincod.html+
\newline
See also W. Brouwer's (less up-to-date) article in 
{\bf Handbook of coding theory}, (ed. Pless et al),
North-Holland, Elsevier, (1998).

\bibitem[E]{E} G. Ewald, {\bf Combinatorial convexity and algebraic geometry}, 
Springer, 1996.

\bibitem[F]{F}
W. Fulton, {\bf Introduction to toric varieties}, PUP, 1993.

\bibitem[GAP]{GAP}
The GAP~Group, {\bf GAP -- Groups, Algorithms, and Programming, 
Version 4.3}; 2000, \verb+(http://www.gap-system.org)+.

\bibitem[GUAVA]{GUAVA}
The GAP package {\bf GUAVA, Version 1.6}, 2002, 
\verb+(http://www.gap-system.org/Info4/deposit.html)+

\bibitem[HCT]{HCT} V. Pless and W. Huffman (eds.),
{\bf Handbook of coding theory}, North-Holland, 1998.

\bibitem[H1]{H} J. Hansen, ``Toric surfaces and error-correcting
codes,'' {\bf Coding theory, cryptography, and related areas},
(ed., Bachmann et al), Springer-Verlag, 1999.

\bibitem[H2]{H2} -----, ``Hirzebruch surfaces and error-correcting
codes,'' preprint, 1999. Available at
\newline
\verb+http://home.imf.au.dk/matjph/+

\bibitem[H3]{H3} -----, ``Toric varieties,
Hirzebruch surfaces and error-correcting
codes,'' preprint, 2002. 

\bibitem[Han1]{Han1} S. Hansen, ``The geometry of Deligne-Lusztig
varieties: higher-dimensional AG codes,'' PhD thesis, 
Univ. Aarhus, 1999 (advisor J. Hansen).

\bibitem[Han2]{Han2} ------, ``Error-correcting codes from
higher-dimensional varieties,'' Finite Fields and their 
Applications \underline{7}(2001)530-552.

\bibitem[Han3]{Han3} ------, ``Algebraic geometry over finite fields - 
with a view towards
applications,'' Master's thesis, Univ. Aarhus, 1996 (advisor J. Hansen).

\bibitem[J]{J} D. Joyner, MAGMA code: \verb+toric.mag+ code,
\newline
\verb+http://web.usna.navy.mil/~wdj/papers/toric.mag+
\newline
GAP code: \verb+toric.g+ code,
\newline
\verb+http://web.usna.navy.mil/~wdj/papers/toric.g+

\bibitem[JV]{JV} ------ and H. Verrill, ``Notes on toric varieties,''
preprint, available at
\newline
\verb+http://web.usna.navy.mil/~wdj/papers.html+

\bibitem[L1]{L1} J. Little, ``A key equation and the computation of error 
values for codes from order domains'',
\newline
\verb+http://www.arXiv.org/abs/math.AC/0303299+

\bibitem[L2]{L2} ------, ``The Ubiquity of Order Domains for 
the Construction of Error Control Codes,'' 
\newline
\verb+http://www.arXiv.org/abs/math.AC/0304292+

\bibitem[MS]{MS} 
F. MacWilliams and N. Sloane, {\bf The theory of error-correcting codes},
North-Holland, 1977.

\bibitem[MAGMA]{MAGMA} 
W. Bosma, J. Cannon, C. Playoust,
``The MAGMA algebra system, I: The user language,''
J. Symb. Comp., \underline{24}(1997)235-265.
\newline
(See also the MAGMA homepage at
\newline
\verb+ http://www.maths.usyd.edu.au:8000/u/magma/+)

\bibitem[S]{S} H. Stichtenoth, {\bf Algebraic function
fields and codes}, Springer-Verlag, 1993.

\bibitem{VT}[VT] M. Tsfasman,
S. Vladut, {\bf Algebraic-geometric codes}, North-Holland, 1998.

\end{thebibliography}
\end{document}